%% file: main.tex
\documentclass{article}
\usepackage{amsthm}
\usepackage{amssymb}
\usepackage{graphicx}
\usepackage{mathtools}
\usepackage{float}
\usepackage{comment}
\usepackage{subfiles}

\newtheorem{theorem}{Theorem}[section]
\newtheorem{lemma}[theorem]{Lemma}
\newtheorem{prop}[theorem]{Proposition}

\newtheorem{cor}[theorem]{Corollary}

\newtheorem{innercustomgeneric}{\customgenericname}
\providecommand{\customgenericname}{}
\newcommand{\newcustomtheorem}[2]{%
  \newenvironment{#1}[1]
  {%
   \renewcommand\customgenericname{#2}%
   \renewcommand\theinnercustomgeneric{##1}%
   \innercustomgeneric
  }
  {\endinnercustomgeneric}
}

\newcustomtheorem{customthm}{Theorem}
\newcustomtheorem{customlemma}{Lemma}
\newcustomtheorem{customprop}{Proposition}
\newcustomtheorem{customcor}{Corollary}

\theoremstyle{remark}

\newcommand{\abs}[1]{\lvert#1\rvert}

\title{Effective Generation of Right-Angled Artin Groups in Mapping Class Groups}
\author{Ian Runnels}
\date{}
\begin{document}

\maketitle

\begin{abstract}
We show that given a collection $X=\{f_1$, \ldots , $f_m\}$ of pure mapping classes on a surface $S$, there is an explicit constant N, depending only on $X$, such that their Nth powers $\{f_1^N$, \ldots , $f_m^N\}$ generate the expected right-angled Artin subgroup of MCG($S$). Moreover, we show that these subgroups are undistorted.
\end{abstract}

\section{Introduction}

\subfile{sections/intro.tex}

\section{Background}

\subfile{sections/background.tex}

\section{Subsurface Projections}

\subfile{sections/subsurfaceprojections.tex}

\section{The Proof of Theorem 4.1}

The goal of this section is to prove

\begin{theorem}
Let $f_1$, \ldots, $f_m$ be an irredundant collection of pure mapping classes supported on connected subsurfaces $S_1$, \ldots , $S_m \subseteq S$. Let
\begin{align*}
 N = \frac{2L + 2 + M_1 + M_2 + 10\delta}{\displaystyle \min_{1 \leq i \leq m} c(S_i)} ,
\end{align*}
where L is as in the Generalized Behrstock inequality, $c(S_i)$ is as in Prop. 3.8, and $M_1$ and $M_2$ are defined below. Then for all $n \geq N$, 
\begin{align*}
   H = \langle f_1^n, \ldots, f_m^n \rangle \cong A(\Gamma),
\end{align*}
where $\Gamma$ is the co-intersection graph of the $S_i$. Moreover, after increasing $N$ in a controlled way, we can guarantee that $H$ is undistorted in $MCG(S)$.
\end{theorem}

We break the proof into two parts, first proving that $H$ is indeed the desired RAAG, then proving that $H$ is undistorted in $MCG(S)$.

\subsection{Theorem 4.1 Part 1: Generation}

\subfile{sections/maingenerated.tex}

\subsection{Theorem 4.1 Part 2: Undistortion}

\subfile{sections/mainundistorted.tex}

\phantom{\cite{farb2011primer} \cite{mousley2018nonexistence}}

\bibliographystyle{plain}
\bibliography{bibliography.bib}

\end{document}

%% file: sections/intro.tex
Let $S = S_{g,p}$ be an orientable surface of genus $g$ with $p$ punctures and satisfying $\chi(S)<0$. By the \textit{mapping class group} of $S$ we mean 
\begin{align*}
    MCG(S) := \pi_0(\textrm{Homeo}^+(S,P)),
\end{align*}
$i.e.$ the group of homotopy classes of orientation-preserving homeomorphisms of $S$ which preserve the set P of punctures. The study of free subgroups of $MCG(S)$ dates back to Klein \cite{klein}, who classically showed that the matrices 
\begin{align*}
    A = \begin{bmatrix} 1 & 2 \\ 0 & 1 \end{bmatrix} \textrm{  and  } \displaystyle B = \begin{bmatrix} 1 & 0 \\ 2 & 1 \end{bmatrix}
\end{align*}
generate a free subgroup of $SL(2,\mathbb{Z})$. Indeed, we may identify $SL(2,\mathbb{Z})$ with $MCG(T^2)$, the mapping class of the torus, and these matrices correspond to the (squares of the) Dehn twists about the standard meridian and longitude curves. It follows from Ivanov's \cite{ivanov1992subgroups} and McCarthy's \cite{mccarthy1985tits} proof of the Tits alternative for $MCG(S)$ that there is in fact an abundance of free subgroups. 

We wish to broaden our view to the larger class of \textit{right-angled Artin groups} (RAAGs). Recall that a RAAG has a presentation determined by a finite simplicial graph $\Gamma$:
\begin{align*}
    A(\Gamma) = \langle v_i \in V(\Gamma) \enskip \vert \enskip [v_i,v_j] = 1 \iff (v_i,v_j) \in E(\Gamma) \rangle.
\end{align*}
Since RAAGs encompass free groups (where $\Gamma$ has no edges), we focus on non-free examples. Regarding these, Koberda showed that they may also be found in abundance. In the statement below, we say that a mapping class is \textit{pure} if it is either pseudo-Anosov or else fixes a multi-curve $C$ component-wise and restricts to either a pseudo-Anosov or the identity on the complementary components $S \backslash C$. We call a mapping class of the latter type a \textit{partial pseudo-Anosov}, and we call the components of $S \backslash C$ where its action is non-trivial its \textit{support}; the support of a pseudo-Anosov is all of $S$.

\begin{theorem}[Koberda, \cite{koberda2012right} Theorem 1.1]
Let $f_1$, \ldots, $f_m$ be an irredundant collection of pure mapping classes supported on connected subsurfaces \linebreak $S_1$, \ldots , $S_m \subseteq S$. There exists some $N \neq 0$ such that for all $n \geq N$, 
\begin{align*}
    \langle f_1^n, \ldots , f_m^n \rangle \cong A(\Gamma),
\end{align*}
where $\Gamma$ is the co-intersection graph of the subsurfaces $\{S_i\}$.
\end{theorem}

Here, \textit{irredundancy} means that no pair of mapping classes have a common power, and the \textit{co-intersection graph} has vertices the subsurfaces $S_i$ and edges between subsurfaces which can be realized disjointly. Koberda's proof goes by playing ping-pong on the space of geodesic laminations on $S$, and it is not clear how the number $N$ depends on the surface $S$ or on the given mapping classes. 

The goal of this paper is to effectivize Koberda's theorem. The constant in the statement of the following theorem below is described in Section 4.
\begin{customthm}{4.1}
Let $f_1$, \ldots, $f_m$ be an irredundant collection of pure mapping classes supported on connected subsurfaces $S_1$, \ldots , $S_m \subseteq S$. There exists a constant $N$, depending explicitly and only on the collection $\{f_i\}$, such that  for all $n \geq N$, 
\begin{align*}
   H = \langle f_1^n, \ldots, f_m^n \rangle \cong A(\Gamma),
\end{align*}
where $\Gamma$ is the co-intersection graph of the $S_i$. Moreover, after increasing $N$ in a controlled way, we can guarantee that $H$ is undistorted in $MCG(S)$.
\end{customthm}

\noindent Computing the constant explicitly in the case that all mapping classes in question are Dehn twists, we have the following corollary.

\begin{customcor}{4.3}
Let $\{t_1$, \ldots, $t_m\}$ be a collection of Dehn twists about distinct curves $\{\beta_1$, \ldots , $\beta_m\}$, and let 
\begin{align*}
     N = 15 + \max_{i,j} i(\beta_i, \beta_j),
\end{align*}
where $i(\cdot, \cdot)$ denotes geometric intersecion number. Then for all $n \geq N$, we have
\begin{align*}
    \langle t_1^n,\ldots,t_m^n \rangle \cong A(\Gamma),
\end{align*}
where $\Gamma$ is the subgraph of $\mathcal{C}(S)$ spanned by \{$\beta_1$, \ldots , $\beta_m$\}.
\end{customcor}

A similar bound has been found by Seo \cite{seo2019powers} using methods from hyperbolic and coarse geometry, and Bass-Serre theory. That these subgroups are undistorted follows from a careful study of the construction of ``admissible" embeddings of RAAGs into mapping class groups due to Clay-Leininger-Mangahas \cite{clay2012geometry}.

It is worth mentioning that if there are more than two mapping classes involved, $N$ \textit{necessarily} depends on the given mapping classes, as the following example illustrates. Let $\alpha_1$ and $\alpha_2$ be two non-trivially intersecting simple closed curves, and consider the Dehn twists
\begin{align*}
    t_1 = t_{\alpha_1}, \enskip t_2 = t_{\alpha_2}, \enskip  \textrm{and} \enskip t_3 = t_1^{2^K} \cdot t_2 \cdot t_1^{-2^K}
\end{align*}
for some $K>0$. Then for no $0< k \leq K$ is $\langle t_1^{2^k}, t_2^{2^k}, t_3^{2^k}  \rangle$ isomorphic to a free group of rank 3, even though the corresponding co-intersection graph is disconnected. 

In Section 2 we establish the relevant notions we will need from coarse geometry, geometric group theory, and the theory of mapping class groups, including a proof of a ping-pong lemma for RAAGs. In Section 3 we recall the theory of subsurface projections due to Masur-Minsky \cite{masur2000geometry}, which we use to build our ping-pong table. The essential result in this section is a modification of the well-known Behrstock inequality \cite{behrstock2006asymptotic}. In Section 4 we carry out the proof of Theorem 4.1.
\section*{Acknowledgements}

The author would like to thank Thomas Koberda, Marissa Loving, and Johanna Mangahas for helpful conversations. The author would also like to thank Rylee Lyman, George Domat, and Jason Behrstock for constructive feedback on the writing of this paper.

%% file: sections/background.tex
\subsection{Quasi-Isometries}

If $(X_1,d_{X_2})$ and $(X_2,d_{X_2})$ are metric spaces and $f:X_1 \to X_2$ is a map, we say $f$ is a \textit{$(A,B)$-quasi-isometric embedding} if there are constants $A \geq 1$ and $B \geq 0$ such that for all $x,y \in X_1$,
\begin{align*}
    \frac{1}{A}d_{X_1}(x,y) - B \leq d_{X_2}(f(x),f(y)) \leq A \cdot d_{X_1}(x,y) + B.
\end{align*}
We also sometimes call such maps ``coarsely Lipschitz". If there is a constant $D>0$ such that any $x_2 \in X_2$ is within $D$ of $f(X_1)$, we further say $f$ is a \textit{quasi-isometry}, and that $X_1$ and $X_2$ are \textit{quasi-isometric}.

Recall that to a group $G$ with generating set $Y$ we may associate the Cayley graph $Cay(G,Y)$, and that equipped with the graph metric $Cay(G,Y)$ is a metric space. Moreover, if $G$ is finitely generated, then any two metrics coming from different finite generating sets $Y$ and $Y'$ yield quasi-isometric Cayley graphs.
We may then put a (left-invariant) metric $d_G$ on $G$, the \textit{word metric}, defined by 
\begin{align*}
    d_{G}(g,h) = d_{Cay(G,Y)}(g,h) = d_{Cay(G,Y)}(1, h^{-1}g).
\end{align*}
Regarding the statement of Theorem 4.1, we say that a subgroup $H < G$ is \textit{undistorted} if the inclusion of $H$ into $G$ is a quasi-isometric embedding with respect to their respective word metrics.

\subsection{Surfaces and their Mapping Class Groups}

Let $S = S_{g,p}$ be a connected, oriented surface of genus $g$ with $p$ punctures (and/or boundary components; the difference here is negligible). The \textit{mapping class group} of $S$, which we denote by $MCG(S)$, is the group of homotopy classes of oritentation-preserving homeomorphisms of $S$ which preserve the set of punctures. We call elements of $MCG(S)$ \textit{mapping classes}. An \textit{essential simple closed curve} is (the homotopy class of) a non-nullhomotopic and non-peripheral simple closed curve, and an \textit{essential subsurface} $S' \subseteq S$ is either a regular neighborhood of an essential simple closed curve, or a component of the complement of a collection of pairwise disjoint essential simple closed curves. As with essential simple closed curves, we consider essential subsurfaces up to homotopy, $i.e.$ up to homotopy of their boundary curves, and in both cases we will not distinguish between an essential simple closed curve/subsurface and any representative. We exclude the pair of pants from our discussion for pathological reasons. 

We frequently study mapping class groups of surfaces by their action on essential simple closed curves and subsurfaces. With respect to this action, we have a ``Jordan canonical form" for mapping classes, due to Thurston \cite{thurston1988geometry}, which has three distinct categories: given $f \in MCG(S)$, $f$ is either finite order, pseudo-Anosov, or reducible, $i.e.$ preserves a multi-curve $C$ in $S$. Though pseudo-Anosov mapping classes are defined geometrically, it follows from this classification that no power of a pseudo-Anosov fixes any essential simple closed curve. For a reducible mapping class $f$, it follows from Birman-Lubotzky-McCarthy \cite{birman1983abelian} that some power $f$ fixes a multi-curve $C$ component-wise, and restricts to a pseudo-Anosov or the identity on each component of $S \backslash C$. We call such a mapping class a \textit{pure reducible} mapping class. In general, by \textit{pure} mapping classes we mean pseudo-Anosovs and pure reducible mapping classes.

\subsection{Right-Angled Artin Groups}

Given a finite simplicial graph $\Gamma$ with vertex set $V(\Gamma)$ and edge set $E(\Gamma)$, the \textit{right-angled Artin group on} $\Gamma$ is the group with presentation
\begin{align*}
    A(\Gamma) := \langle v_i \in V(\Gamma) \enskip \vert \enskip [v_i,v_j] = 1 \iff (v_i,v_j) \in E(\Gamma) \rangle.
\end{align*}
We call the $v_i$ the \textit{vertex generators} of $A(\Gamma)$. The standard examples of such groups are free groups (where $\Gamma$ has no edges), free abelian groups (where $\Gamma$ is a complete graph), and free and direct products of such groups (corresponding to disjoint union of graphs and join of graphs, respectively). Despite their simple presentations, these groups can be quite complicated. That said, they are still universal enough that we can study their actions on spaces quite easily. The following is a modification of the ping-pong lemma for RAAGs found in \cite{koberda2012right}, and is the main tool in proving Theorem 4.1.

\begin{lemma}[Ping-Pong]
Let $A(\Gamma)$ be a right-angled Artin group on a graph $\Gamma$ which is anti-connected, $i.e.$ the complement graph $\Gamma^c$ is connected. Suppose $A(\Gamma)$ is acting on a set $X$ such that there exist non-empty subsets $X_i' \subset X_i \subset X$ for each vertex generator $v_i$ satsifying
\begin{enumerate}
    \item If $X_i \cap X_j \neq \varnothing$, then there exists $x_i \in X_i$ which does not belong to $X_j$, and vice versa
    \item If $u$ is a word not containing a power of $v_j$, wherein every vertex generator commutes with $v_j$, then $u(X_j') \subset X_j$ 
    \item If $v_i$ and $v_j$ do not commute, then $X_i$ and $X_j$ are disjoint and $v_i^m(X_j) \subset X_i'$ for all $m \neq 0$
\end{enumerate}
\noindent Then the $A(\Gamma)$ action on X is faithful.
\end{lemma}

\begin{proof}
If $\Gamma$ is not anti-connected, then $A(\Gamma)$ splits as a direct product, and we can play ping-pong on each factor. We begin by putting $w$ in a normal form called \textit{central form}, due to M. Kapovich ($cf.$ \cite{kapovich2012raags},\cite{koberda2012right}). Using only commutations, we may write $w$ as
\begin{align*}
    w = u_k \cdot v_{i_k}^{r_k} \cdot u_{k-1} \cdot v_{i_{k-1}}^{r_{k-1}} \cdots u_1 \cdot v_{i_1}^{r_1}
\end{align*}
where
\begin{itemize}
    \item each $u_j$ is a word in the vertex generators of $A(\Gamma)$, such that every generator appearing in $u_j$ commutes with every other generator appearing in $u_j$
    \item $v_{i_j}$ does not commute with $v_{i_{j+1}}$
    \item $v_{i_j}$ commutes with every generator appearing in $u_j$
\end{itemize}
\noindent We will induct on the central-word-length of reduced words in $A(\Gamma)$. For the base case, we have $w = u_1 v_{i_1}$. Since we assumed $\Gamma$ was anti-connected, there is some generator $v_j$ which does not commute with $v_{i_1}$. Choosing $x_j \in X_j$ and applying (3) we have $v_{i_1}^{r_1} \cdot x_j \in X_{i_1}'$. Applying (2), we then have $u_1 \cdot v_{i_1}^{r_1} \cdot x_j \in X_{i_1}$. Again by (3), since $X_{i_1} \cap X_j = \varnothing$, we see that $w \cdot x_j \neq x_j$ and we are done. Inductively, we have 
\begin{align*}
w = u_k \cdot v_{i_k}^{r_k} \cdots u_2 \cdot v_{i_2}^{r_2} \cdot u_1 \cdot v_{i_1}^{r_1}
\end{align*}
By either (1) or (3), we can choose $x_{i_2} \in X_{i_2}$ which does not belong to $X_{i_k}$; note that since $v_{i_2}$ and $v_{i_1}$ don't commute by construction, $x_{i_2}$ also does not belong to $X_{i_1}$. Repeatedly applying the argument above to this element, we have by induction that $w \cdot x_{i_2} \in X_{i_k}$ , so in particular $w \cdot x_{i_2} \neq x_{i_2}$.
\end{proof}

%% file: sections/subsurfaceprojections.tex
\noindent Our ping-pong sets will be given in terms of \textit{subsurface projections} of simple closed curves to subcomplexes of the curve complex of $S$, as originally defined by Masur-Minsky \cite{masur2000geometry}. Recall that the \textit{curve complex} of $S$, denoted $\mathcal{C}(S)$, is the simplicial complex whose vertices are homotopy classes of simple closed curves, and whose simplices correspond to collections of curves which can be realized (pairwise) disjointly. We equip $\mathcal{C}(S)$ with the graph metric. A celebrated theorem of Masur-Minsky \cite{masur1999geometry} says that with this metric, $\mathcal{C}(S)$ is $\delta$-hyperbolic. Moreover, Aougab \cite{aougab2013uniform}, Bowditch \cite{bowditch2014uniform}, and Clay-Rafi-Schleimer \cite{clay2015uniform} showed that $\delta$ can be made independent of $S$, and Hensel-Przytycki-Webb \cite{hensel20151} showed that $\delta$ can be taken to be 17.

\subsection{Constructing Subsurface Projections}

Given a $\pi_1$-injective, non-annular subsurface $S' \subset S$ and a simple closed curve $\gamma$, we define the projection $\pi_{S'}(\gamma)$ of $\gamma$ to $S'$ as follows: if $\gamma$ is disjoint from $S'$, then $\pi_{S'}(\gamma) = \varnothing$, and if $\gamma$ is entirely contained is $S'$ then $\pi_{S'}(\gamma) = \gamma$. Otherwise, $\gamma$ non-trivially intersects the boundary curves of $S'$, and we define $\pi_{S'}(\gamma)$ to be the set of essential curves obtained by taking the boundary of a neighborhood of $(\gamma \cap S') \cup \partial S'$, see Figure 1.

\begin{figure}[H]
\centering
  \includegraphics[width=.85\linewidth]{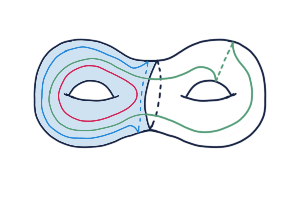}
  \caption{The projection of the green curve to the left one-holed torus consists of the red and blue curves.}
  \label{Fig. 1}
\end{figure}

There is always at least one such curve, and we think of these projections as living in $\mathcal{C}(S')$, as a subcomplex of $\mathcal{C}(S)$. Given two curves $\beta$ and $\gamma$ with non-trivial projection to $S'$, we define their \textit{projection distance} $d_{S'}(\beta, \gamma)$ to be 
\begin{align*}
    d_{S'}(\beta, \gamma) = \inf\{ d_{S'}(b,c) \enskip \vert \enskip b \in \pi_{S'}(\beta), c \in \pi_{S'}(\gamma)\}.
\end{align*}
We note that this definition differs from the original, which uses Hausdorff distance instead.

The projection of a curve $\gamma$ to an annulus about another curve $\beta$ is defined slightly differently: we first fix a hyperbolic metric on $S$ and take geodesic representatives of $\gamma$ and $\beta$. Consider the cover $S_{\beta}$ of $S$ corresponding to the cyclic subgroup of $\pi_1(S)$ generated by $\beta$, which we think of as an infinite flaring annulus with core curve $\widetilde{\beta}$, compactified with its Gromov boundary. We define the projection $\pi_{\beta}(\gamma)$ to be the collection of lifts $c$ of $\gamma$ to the cover $S_{\beta}$ which cross the core curve $\widetilde{\beta}$, connecting the two boundary components, see Figure 2. We can assemble the set of all homotopy (rel. boundary) classes of arcs in $S_{\beta}$ into a simplicial complex $\mathcal{A}(\beta)$, the \textit{arc complex} of $\beta$, with edges representing arcs with disjoint interiors, and equipped with the graph metric. Note that given any arc $\gamma$ in $\mathcal{A}(\beta)$, all of its parallel translates around $S_{\beta}$ (of which there are uncountably many) are by definition contained in the 1-neighborhood of $\gamma$, and it is easy to show that $\mathcal{A}(\beta)$ is quasi-isometric to $\mathbb{Z}$. Given another curve $\delta$ crossing $\beta$, we define the projection distance $d_{\beta}(\gamma, \delta)$ in the same way as above. Though we chose a hyperbolic metric, it is not hard to see that $d_{\beta}(\delta,\gamma) \leq i(\delta,\gamma)$.

\begin{figure}[H]
\centering
  \includegraphics[width=.85\linewidth]{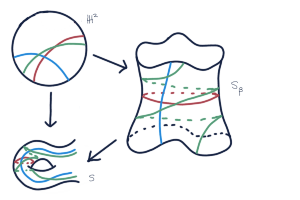}
  \caption{The annular cover of S corresponding to the red curve, with lifts of the blue and green curves.}
  \label{Fig. 2}
\end{figure}

One useful fact (\cite{masur2000geometry},\cite{behrstock2006asymptotic}) about these projections is that they are coarsely Lipschitz. Another concerns the effect of mapping classes supported away from the subsurface in question.

\begin{lemma}
Let $f \in$ MCG(S) be a pure mapping class supported on a subsurface $S_j$ disjoint from $S_i$. If $S_i$ is an annulus about a simple closed curve $\beta$, we also require that $\partial S_j$ does not contain $\beta$. Then, if $S_i$ is not an annulus, we have $d_{S_i}(\delta, f(\gamma)) = d_{S_i}(\delta,\gamma)$. If $S_i$ is an annulus about a curve $\beta$ then $d_{\beta}(\delta, f(\gamma))$ differs from $d_{\beta}(\delta,\gamma)$ by at most 2.
\end{lemma}

\begin{proof}
For a non-annular subsurface $S_i$, observe that if $f$ is supported disjointly from $S_i$, then $f$ is either the identity on $S_i$ or it acts as a power of a Dehn twist along the boundary $\partial S_i$. If $\gamma$ is a curve disjoint from $\partial S_i$, then $f(\gamma)$ is as well, so the projection is clearly unchanged. Otherwise, the intersection of $f(\gamma)$ with $S_i$ differs from that of $\gamma$ with $S_i$ only by twisting about the boundary $\partial S_i$, which doesn't affect the surgery construction of the projection. Thus, $\pi_{S_i}(f\cdot \gamma) = \pi_{S_i}(\gamma)$, and so projection distances remain unchanged. 

For projections to annuli, this is most easily seen by lifting the action of $f$ to $S_{\beta}$. We will give the argument in the case that $f$ is a Dehn twist about another simple closed curve $\eta$, but the general case is the same. Since $\eta$ is disjoint from $\beta$, no lift of $\eta$ crosses the core curve of $S_\beta$, so any such lift must connect a boundary component of $S_{\beta}$ to itself. Moreover, since we assumed $\eta$ was simple, each lift has distinct endpoints in the boundary and no two lifts cross each other. Hence, these lifts cannot traverse a full circuit about the boundary of $S_{\beta}$. Then the action of $f$  on lifts of $\gamma$ is by shifting their endpoints along the top and bottom boundary components of $S_{\beta}$. Since this shifting cannot move endpoints more than a full revolution around each boundary component, the difference between $d_{\beta}(\delta,\gamma)$ and $d_{\beta}(\delta, f(\gamma))$ is at most 2.
\end{proof}

It is worth noting that it is precisely because of this ``coarse invariance" of annular projection distance that, in the statement of Lemma 2.1 above, we required the existence of ``coarsely preserved" subsets $X_i' \subset X_i \subset X$. Since non-annular projections are preserved on the nose, when considering other pure mapping classes these subsets are unnecessary.

\subsection{The Masur-Minsky Distance Formula}

To show that the RAAGs we generate are undistorted in $MCG(S)$, we will need a way of relating word length in $MCG(S)$ to subsurface projection distances. This relationship is captured by the ``Masur-Minsky distance formula" of \cite{masur2000geometry}. Before stating it, we establish notation. A (complete clean) \textit{marking} $\mu$ on $S$ consists of a pants decomposition $\{\beta_1,\ldots,\beta_{3g-3+p}\}$, called the \textit{base} of $\mu$, together with a \textit{transversal} for each $\beta_i$ satisfying certain properties which are unnecessary for the discussion at hand. Masur-Minsky build a simplicial complex $\widetilde{\mathcal{M}}(S)$, called the \textit{marking complex} of $S$, whose vertex set is set of all markings and whose edges correspond to certain \textit{elementary moves} on markings. Equipped with the graph metric, the complex $\widetilde{\mathcal{M}}(S)$ is locally finite and admits an action of $MCG(S)$ by isometries, so that $MCG(S)$ and $\widetilde{\mathcal{M}}(S)$ are quasi-isometric.
We define the projection of a marking $\mu$ to an essential non-annular subsurface $S'$ to be $\pi_{S'}(base(\mu))$. We define the projection of $\mu$ to an essential annulus to be either $\pi_{S'}(base(\mu))$ if the core curve of the annulus is not in base($\mu$), and the projection of the transversal otherwise. We can now state
\begin{theorem}[Masur-Minsky, \cite{masur2000geometry}]
There exists $K_0 = K_0(S) > 0$ such that for all $K \geq K_0$, there exist constants $A \geq 1$ and $B \geq 0$ with the property that for all pairs of markings $\mu, \mu' \in \widetilde{\mathcal{M}}(S)$ we have
\begin{align*}
    \frac{1}{A}\sum_{S' \subseteq S} [[d_{S'}(\mu,\mu')]]_K - B \leq d_{\widetilde{\mathcal{M}}(S)}(\mu,\mu') \leq A\sum_{S' \subseteq S} [[d_{S'}(\mu,\mu')]]_K +B,
\end{align*}
where $[[x]]_K = x$ if $x \geq K$ and is 0 otherwise.
\end{theorem}

In particular, we can approximate the word length of a mapping class $f$ by looking at the subsurface projections distances between $\mu$ and $f(\mu)$.

\subsection{A Generalized Behrstock Inequality}

\noindent A key idea in the proof of Theorem 4.1 is a slight generalization of the following inequality due to Behrstock \cite{behrstock2006asymptotic}. A constructive proof can be found in \cite{mangahas2013recipe}.

\begin{lemma}[Behrstock Inequality]
Let $S_i$, $S_j$, and $S_k$ be three pairwise intersecting essential subsurfaces or simple closed curves. Then
\begin{align*}
    d_{S_i}(\partial S_j, \partial S_k) \geq 10 \implies d_{S_j}(\partial S_i,\partial S_k) \leq 4.
\end{align*}
If $S_i$ (or $S_j$ or $S_k$) is an annulus, we replace $\partial S_i$ with the core curve $\alpha_i$. If all three are annuli, we may further replace 4 with 3.
\end{lemma}

The generalization we will use concerns not only subsurface projections, but nearest point projections to geodesics in $\mathcal{C}(S)$ and its subcomplexes. For the proof, we will need the following two results. The first is a straightforward computation in $\delta$-hyperbolic geometry.

\begin{prop}
Let $\alpha \subset \mathcal{C}(S)$ be a geodesic, and $x,y \in \mathcal{C}(S)$. Then
\begin{align*}
    d(\pi_{\alpha}(x),\pi_{\alpha}(y)) \geq 8 \delta + 2 \implies  d(\pi_{\alpha}(x),\pi_{\alpha}(y)) \leq d(x,y) + 24\delta,
\end{align*}
where $\pi_{\alpha}$ is a (coarse) nearest point projection map.
\end{prop}

The second is a theorem of Masur-Minsky \cite{masur2000geometry}, known as the Bounded Geodesic Image Theorem (BGIT). The uniform, effective statement below is due to Webb \cite{webb2015uniform}.

\begin{theorem}[BGIT]
If $S' \subset S$ is a subsurface and $\alpha$ is a geodesic in $\mathcal{C}(S)$ with the property that $\pi_{S'}(z) \neq \varnothing$ for all $z \in \alpha$, then
\begin{align*}
    diam_{S'}(\pi_{S'}(\alpha)) < 100.
\end{align*}
\end{theorem}

The lemma below was originally observed by Sun \cite{sunthesis}. We provide an original proof for the reader's convenience. In the statement, by $d_{\alpha_i}(.,.)$, we mean $d(\pi_{\alpha_i} \circ \pi_{S_i}(.),\pi_{\alpha_i} \circ \pi_{S_i}(.))$.

\begin{lemma}[Generalized Behrstock Inequality]
Let $\gamma$ be a non-peripheral simple closed curve and let $\alpha_1$, $\alpha_2$ be either non-peripheral simple closed curves or geodesics in $\mathcal{C}(S_1) \subseteq \mathcal{C}(S)$, $\mathcal{C}(S_2) \subseteq \mathcal{C}(S)$, respectively. If $\alpha_1$ and $\alpha_2$ are both geodesics, we assume that they do not share endpoints in the boundary of $\mathcal{C}(S)$. Then there is an $L = L(\delta) > 0$ such that 
\begin{align*}
    \min\{d_{\alpha_1}(\gamma, \alpha_2) , \enskip d_{\alpha_2}(\gamma, \alpha_1)\} \leq L.
\end{align*}
\end{lemma}

\begin{proof}
We break the proof into cases depending on whether $\alpha_1$ and $\alpha_2$ are simple closed curves and/or geodesics, and depending on the topological configuration of $S_1$ and $S_2$. If $\gamma$ is disjoint from either $S_1$ or $S_2$ or if $S_1$ and $S_2$ are disjoint, there is nothing to prove, so we assume these intersections are all non-empty. The requirement that two geodesics not share endpoints ensures that they only fellow travel for a finite amount of time, which by hyperbolicity implies that each has bounded projection to the other. We note that in the arguments below there is repeated implicit use of Proposition 3.4 and the fact that projection distance is definied as an infimum.

\noindent \underline{\textbf{Case 1}}: $\alpha_1$, $\alpha_2$ are simple closed curves.

\vspace{\baselineskip}

In this case, we can just use the Behrstock Inequality.

\vspace{\baselineskip}

\noindent \underline{\textbf{Case 2}}: $\alpha_1$ is a simple closed curve and $\alpha_2$ is a geodesic in $\mathcal{C}(S_2)$.

\vspace{\baselineskip}

Suppose $d_{\alpha_1}(\gamma, \alpha_2) > 100$, $i.e.$ $d_{\alpha_1}(\gamma, z) > 100$ for all $z \in \alpha_2$ with non-trivial projection to $\alpha_1$. By the contrapositive of BGIT, this implies that any geodesic between $\gamma$ and any $z \in \alpha_2$ passes through the 1-neighborhood of $\alpha_1$. In particular, on any geodesic between $\gamma$ and $\pi_{\alpha_2}(\gamma)$, there is a vertex $x$ with $d(x,\alpha_1) = 1$. But the projection of $x$ to $\alpha_2$ overlaps with that of $\gamma$, so we have
\begin{align*}
    d_{\alpha_2}(\gamma, \alpha_1) & \leq d_{\alpha_2}(\gamma, x) + d_{\alpha_2}(x,\alpha_1) \\
    &\leq 0 + (1 + 24\delta).
\end{align*}

\noindent\underline{\textbf{Case 3}}: $\alpha_1$ and $\alpha_2$ are both geodesics. 

\vspace{\baselineskip}

We first consider the case that $S_1 = S_2$. Suppose $d_{\alpha_2}(\gamma,\alpha_1) > 28\delta$, so that $d(\gamma,\alpha_1) > 4\delta$. By hyperbolicity, a geodesic in $\mathcal{C}(S_2)$ connecting $\gamma$ and $\pi_{\alpha_1}(\gamma)$ passes within 2$\delta$ of $\alpha_2$, $i.e.$ there is a vertex $z$ on this geodesic and a vertex $y$ on $\alpha_2$ with $d(z,y) \leq 2\delta$. Then
\begin{align*}
    d_{\alpha_1}(\gamma,\alpha_2) &\leq d_{\alpha_1}(\gamma,z) + d_{\alpha_1}(z,y) + d_{\alpha_1}(y,\alpha_2) \\
    &\leq 0 + (2\delta + 24\delta) + 0.
\end{align*}
If $S_1 \subset S_2$, then since subsurface projections are coarsely Lipschitz, we can perform exactly the same argument provided we assume $d_{\alpha_1}(\gamma,\alpha_2) > 28\delta$. Lastly, we suppose $\partial S_1$ and $\partial S_2$ intersect. If $\gamma$ lies entirely within $S_1$, then
\begin{align*}
    d_{\alpha_2}(\gamma,\alpha_1) \leq d_{\alpha_2}(\gamma, \partial S_1) + d_{\alpha_2}(\partial S_1,\alpha_1).
\end{align*}
As $\partial S_1$ is disjoint from both $\gamma$ and $\alpha_1$, each term on the right is bounded (in terms of $\delta$ only). If $\gamma$ intersects both $\partial S_1$ and $\partial S_2$, and $d_{\alpha_1}(\gamma,\alpha_2) > 11 + 48\delta$, then
\begin{align*}
    d_{S_1}(\gamma,\partial S_2) &\geq d_{\alpha_1}(\gamma,\partial S_2) - 24\delta \\
    &\geq d_{\alpha_1}(\gamma,\alpha_2) - d_{\alpha_1}(\alpha_2, \partial S_2) - 24\delta \\
    &\geq (11 + 48\delta) - (1 + 24\delta) - 24\delta \\
    &\geq 10.
\end{align*}
Hence, by the Behrstock Inequality, $d_{S_2}(\gamma,\partial S_1) \leq 4$, and so
\begin{align*}
    d_{\alpha_2}(\gamma,\alpha_1) &\leq d_{\alpha_2}(\gamma, \partial S_1) + d_{\alpha_2}(\partial S_1,\alpha_1) \\
    &\leq (d_{S_2}(\gamma,\partial S_1) + 24\delta) + (1 + 24\delta) \\
    &\leq 5 + 48\delta.
\end{align*}

\noindent It follows in each case that taking $L$ larger than all of the constants appearing in each case, we get the desired conclusion.

\end{proof}

An immediate consequence of the Lemma 3.6 is that we obtain a partial ordering on the set of geodesics in $\mathcal{C}(S)$ and its subcomplexes, in analogy with the partial order on subsurfaces in \cite{behrstock2012geometry}, afforded by the Behrstock inequality. Indeed, for $K \geq 2L$ let $\Omega(K,\mu,\mu')$ denote the set of such geodesics with $d_{\alpha}(\mu,\mu') \geq K$, and let $\alpha_1$ and $\alpha_2 \in \Omega(K,\mu,\mu')$. We will say $\alpha_1 \prec \alpha_2$ if $d_{\alpha_1}(\mu, \alpha_2) \geq L$. We then have the following characterization of this ordering.

\begin{prop}
Let $\alpha_1, \alpha_2 \in \Omega(K,\mu,\mu')$. Then $\alpha_1$ and $\alpha_2$ are ordered and the following are equivalent
\begin{enumerate}
    \item $\alpha_1 \prec \alpha_2$;
    \item $d_{\alpha_1}(\mu,\alpha_2) \geq L$;
    \item $d_{\alpha_1}(\mu',\alpha_2) < L$;
    \item $d_{\alpha_2}(\mu',\alpha_1) \geq L$;
    \item $d_{\alpha_2}(\mu, \alpha_1) < L$.
\end{enumerate}
\end{prop}

The proof is the same as that of \cite{clay2012geometry}, Prop. 3.6, replacing the Behrstock inequality with the Generalized Behrstock inequality where appropriate.

\subsection{The Action on the Curve Complex}

The following are a pair of propositions of Masur-Minsky \cite{masur1999geometry}\cite{masur2000geometry} concerning the action on the curve complex of a pseudo-Anosov mapping class. The first tells us that pseudo-Anosovs act on $\mathcal{C}(S)$ like hyperbolic isometries.

\begin{prop}[Masur-Minsky \cite{masur1999geometry}, Prop. 3.6]
There exists a constant $c = c(S)>0$ such that, for any pseudo-Anosov mapping class $f \in MCG(S)$, any simple closed curve $\gamma$, and any $n \in \mathbb{Z} \backslash \{0\}$, we have
\begin{align*}
    d_S(f^n(\gamma),\gamma) \geq c \abs{n} .
\end{align*}

\end{prop}
Masur-Minsky proved the above for surfaces satisfying $3g - 3 + p >4$. For the so-called ``sporadic" cases, namely $S_{1,1}$ and $S_{0,4}$, we redefine the curve complex in such a way that we obtain the Farey graph, where it is noted by Mangahas \cite{mangahas2010uniform} that the same result follows by considering the action of hyperbolic isometries on the Farey graph embedded in $\mathbb{H}^2$. It is easy to show that Proposition 3.6 implies that for any simple closed curve $\gamma$ and any pseudo-Anosov $f$, the bi-infinite sequence of curves $\{f^n(\gamma) \vert n \in \mathbb{Z}\}$ is an $f$-invariant quasi-geodesic. By restricting a pure mapping class to a pseudo-Anosov component $S' \subset S$, we obtain such a lower bound for the action of $f$ on the curve complex of $S'$ as a subcomplex of that of $S$, and for a power of a Dehn twist, the quantity $c$ can be taken to be 1. 

As noted above, any pseudo-Anosov $f$ preserves many \textit{quasi-geodesics} in $\mathcal{C}(S)$. However, the Generalized Behrstock Inequality was stated in terms of projections to \textit{geodesics}. In order to apply the Generalized Behrstock Inequality, we will need the following

\begin{prop}[Masur-Minsky \cite{masur2000geometry}, Prop. 7.6]
Let $f \in MCG(S)$ be pseudo-Anosov. There exists a bi-infinite geodesic $\beta$ in $\mathcal{C}(S)$ such that for all $j$, $\beta$ and $f^j(\beta)$ are $2\delta$ fellow travelers. 
\end{prop}

The geodesic $\beta$ and its $f$-translates are referred to as a \textit{quasi-axis} for $f$. Applying Proposition 3.6 to the action of $f$ on its quasi-axis, we have

\begin{lemma}[Masur-Minsky \cite{masur2000geometry}, Lemma 7.7] Given $A > 0$, let $N$ be the smallest integer such that $c(S) \cdot N > A + 10\delta$, where $c(S)$ is the constant from Proposition 3.6. Then for all $n \geq N$, 
\begin{align*}
    d(\pi(x),\pi(f^n(x))) \geq A .
\end{align*}
where $\pi$ denotes projection to the quasi-axis of $f$.
\end{lemma}

%% file: sections/maingenerated.tex
We first show that the group generated by $f_1^N,\ldots,f_m^N$ is the expected RAAG.

\begin{proof}
Let $\{f_1$, \ldots, $f_m\} \in$ MCG(S) be an irredundant collection of pure mapping classes with supporting subsurfaces $\{S_1$, \ldots, $S_m\}$, and for each $i$ let $\alpha_i$ be a geodesic for $f_i$ in $\mathcal{C}(S_i) \subseteq \mathcal{C}(S)$ as in Prop. 3.9, or the core curve of $S_i$ if $S_i$ is an annulus. Without loss of generality, we assume that the co-intersection graph of the $S_i$ is anti-connected, so that for each $f_i$ there is some $f_j$ which does not commute with it. We will explicitly construct a constant $N$ and a group action so that $\{f_1^N$, \ldots, $f_m^N\}$ satisfy the criteria for ping-pong. To this end, set
\begin{align*}
    X = \{\gamma \enskip \vert \enskip \gamma \textrm{ an essential simple closed curve in }S\} ,
\end{align*}
and for each $1 \leq i \leq m$, set
\begin{align*}
    X_i &= \{\gamma  \enskip \vert \enskip d_{\alpha_i}(\gamma,\alpha_j) > L \enskip \forall \enskip j \textrm{ such that } S_j \cap S_i \neq \emptyset \} ,\\
    X_i' &= \{\gamma  \enskip \vert \enskip d_{S_i}(\gamma,\partial S_j) > L + 2 \enskip\forall \enskip j \textrm{ such that } S_j \cap S_i \neq \emptyset \} ,
\end{align*}
where $L$ is as in the Generalized Behrstock Inequality. Observe that by the Generalized Behrstock Inequality, if $S_i$ and $S_j$ cannot be realized disjointly, then their corresponding sets $X_i$ and $X_j$ are disjoint. Moreover, since we assumed the mapping classes were irredundant, no pair $\alpha_i,\alpha_j $ can share endpoints in the boundary of $\mathcal{C}(S)$.

Let $w \in \langle f_1^N,\ldots,f_m^N \rangle$. We begin by putting $w$ in central form as in Lemma 2.1. Using only commutations, we may write $w$ as
\begin{align*}
    w = u_k \cdot g_k \cdot u_{k-1} \cdot g_{k-1} \cdots u_1 \cdot g_1,
\end{align*}
where each $g_j$ represents some power of some $f_i^N$, and each $u_j$ is a word in $\langle f_1^N,\ldots,f_m^N \rangle$. We (possibly) make one further modification to this representative of $w$. For each $g_j$ which is a power of a Dehn twist, if a generator appearing in the corresponding $u_j$ is supported on a subsurface containing the twisting curve as a boundary component, we compose $g_j$ with this mapping class to obtain a single mapping class with the same support, and rearrange $w$ to be in central form once again. Hence, in this final representative for $w$, if any $g_j$ is a Dehn twist, then no mapping class appearing in $u_j$ can alter projections to the corresponding annulus by much (Lemma 3.1).

We may now play ping-pong. Write $w$ in central form as above. Up to relabelling, we assume $g_1 = f_1^{Nr_1}$, $g_2 = f_2^{Nr_2}$, and $g_k = f_j^{Nr_j}$ for some $j$. Choose $\gamma \in X_2 \backslash(X_2 \cap X_j)$; either $g_2$ and $g_k$ don't commute, so their corresponding sets $X_2$ and $X_j$ are disjoint, or they commute and their supports are disjoint, and we can choose some curve $\gamma$ which intersects the support of $g_2$ but not that of $g_k$. If $g_k$ is also a power of $f_2^N$, conjugate $w$ by $g_k$, choose $\gamma \in X_{1}$, and run the same argument below. Since $g_1$ and $g_2$ don't commute, their corresponding sets $X_1$ and $X_2$ are disjoint. In particular, since $\gamma \in X_2$, it satisfies
\begin{align*}
    d_{\alpha_1}(\gamma, \alpha_2) < L.
\end{align*}
Applying the triangle inequality, we then have
\begin{align*}
    d_{\alpha_1}(\gamma, \alpha_{\ell}) &\leq d_{\alpha_1}(\gamma, \alpha_2) + d_{\alpha_1}(\alpha_2, \alpha_{\ell}) \\
    &\leq L + M_1,
\end{align*}
where $M_1 = \displaystyle \max_{1 \leq i,\ell,r \leq m} d_{\alpha_i}(\alpha_{\ell}, \alpha_r)$. Applying the triangle inequality again, 
\begin{align*}
    d_{\alpha_1}(f_1^N(\gamma), \alpha_{\ell}) &\geq d_{\alpha_1}(\gamma, f_1^N(\gamma)) - d_{\alpha_1}(\gamma, \alpha_{\ell}) \\
    &\geq d_{\alpha_1}(\gamma, f_1^N(\gamma)) -L -M_1.
\end{align*}
Hence, if $d_{\alpha_1}(\gamma, f_1^N(\gamma)) \geq 2L + 2 + M_1 + M_2$, where
\begin{align*}
    M_2 = \max_{1 \leq i,j \leq m} diam\{d_{\alpha_i}(\alpha_j)\},
\end{align*}
we will have $f_1^N(\gamma) \in X_1'$. Invoking Lemma 3.10, we set
\begin{align*}
    N = \frac{2L + 2 + M_1 + M_2 + 10\delta}{\displaystyle \min_{1 \leq i \leq m} c(S_i)}.
\end{align*}
Thus, $g_1(\gamma) \in X_1'$, and by Lemma 3.1, $u_1 \cdot g_1(\gamma) \in X_1$. Inductively, we then have $w(\gamma) \in X_j$, so $w$ acts non-trivially on $X$ and we are done.
 
\end{proof}

We have two immediate corollaries, which follow from chasing the appropriate constants through the above argument.

\begin{cor}
For any pair of pseudo-Anosov and/or partial pseudo-Anosov mapping classes $f$ and $g$, if $\langle f,g \rangle$ is not virtually abelian, then there exists $N = N(S)$ depending only on $S$ such that for all $n \geq N$
\begin{align*}
    \langle f^n, g^n \rangle \cong F_2 .
\end{align*}
\end{cor}

The second corollary is simply the restriction of Theorem 1.3 to the case where all mapping classes are Dehn twists:

\begin{cor}
Let $\{t_1$, \ldots, $t_m\}$ be a collection of Dehn twists about distinct curves $\{\beta_1$, \ldots , $\beta_m\}$, and let 
\begin{align*}
     N = 15 + \max_{i,j} i(\beta_i, \beta_j).
\end{align*}
Then for all $n \geq N$, we have
\begin{align*}
    \langle t_1^n,\ldots,t_m^n \rangle \cong A(\Gamma),
\end{align*}
where $\Gamma$ is the subgraph of $\mathcal{C}(S)$ spanned by \{$\beta_1$, \ldots , $\beta_m$\}.
\end{cor}

This should be compared to the main theorem of \cite{seo2019powers}, where a similar (quadratic) bound was computed. As an easy example, we have the following.
\begin{cor}
The $16^{th}$ powers of the Humphries generators for $MCG(S)$, whose underlying curves are pictured below, generate a RAAG.
\end{cor}

\begin{figure}[H]
\centering
  \includegraphics[width=.75\linewidth, height=.4\linewidth]{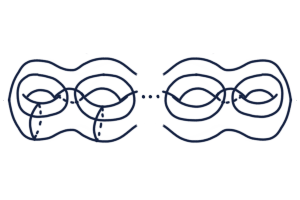}
  \caption{The curves whose Dehn twists are the Humphries generators of $MCG(S)$}.
  \label{Fig. 3}
\end{figure}

%% file: sections/mainundistorted.tex
We now show that the subgroups generated in the previous section are undistorted in $MCG(S)$, after raising the power $N$ by a controlled amount. To do this, we will use the following theorem from \cite{clay2012geometry}. Though we use projections to geodesics instead of just subsurface projections, the proof is identical (and somewhat technical, so we refer the reader to the source).

\begin{theorem}
Let $H$ be as above, $\mu \in \widetilde{\mathcal{M}}(S)$ be a marking on $S$, and let
\begin{align*}
    N = \frac{2L + 2 + M_1 + M_2 + 10\delta + 2M_3 + K_0}{\displaystyle \min_{1 \leq i \leq m} c(S_i)},
\end{align*}
where
\begin{align*}
    M_3 = \max_{1 \leq i,j \leq m} d_{\alpha_i}(\mu,\alpha_j),
\end{align*}
and where $K_0$ is as in the Masur-Minsky distance formula. Set $K$ to be the numerator of N. Let $w = u_k \cdot g_k \cdots u_1 \cdot g_1 \in H$. Then 
\begin{align*}
    d_{u_k \cdot g_k \cdots u_{i-1} \cdot g_{i-1}\alpha_j}(\mu, w \cdot \mu) \geq K\abs{e_i} ,
\end{align*}
where $g_i = (f_j^N)^{e_i}$.
\end{theorem}

The proof of (un)distortion also closely follows \cite{clay2012geometry}, with one small modification.
\begin{proof}
Via the quasi-isometry between $MCG(S)$ and $\widetilde{\mathcal{M}}(S)$, it suffices to show that there are constants $A \geq 1$ and $B \geq 0$ such that for all $w \in H$
\begin{align*}
    \frac{1}{A}d_{\widetilde{\mathcal{M}}(S)}(\mu,w \cdot \mu) - B \leq d_{H}(1,w) \leq A \cdot d_{\widetilde{\mathcal{M}}(S)}(\mu,w \cdot \mu) + B.
\end{align*}
For any group $G$ acting by isometries on a metric space $(X,d_X)$, we always have 
\begin{align*}
    d_X(x, g \cdot x) \leq A \cdot d_G(1,g),
\end{align*}
where $A \geq \displaystyle\max d_{X}(x, s_i \cdot x)$, and $s_i$ is a generator for $G$. Hence, we need only to find $A$ and $B$ so that for all $w \in H$
\begin{align*}
    d_H(1, w) \leq A \cdot d_{\widetilde{\mathcal{M}}(S)}(\mu,w \cdot \mu) + B.
\end{align*}
Let $H$, $N$, and $K$ be as above, and let $w = u_k \cdot (f_{\ell_k}^N)^{e_k} \cdots u_1 \cdot (f_{\ell_1}^N)^{e_1}$ and $w_i = u_k \cdot (f_{\ell_k}^N)^{e_k} \cdots u_i \cdot (f_{\ell_i}^N)^{e_i} $. Then
\begin{align*}
    d_{H}(1,w) &= \sum_{i=1}^k \abs{e_i} \\
    &\leq \sum_{i=1}^k K\abs{e_i} \\
    &\leq \sum_{i=1}^k d_{w_{i-1} \cdot \alpha_i}(\mu, w \cdot \mu). \\
\end{align*}
By Proposition 3.4, each term in the last sum satisfies
\begin{align*}
    d_{w_{i-1} \cdot \alpha_i}(\mu, w \cdot \mu) \leq d_{w_{i-1} \cdot S_i}(\mu, w \cdot \mu) + 24\delta.
\end{align*}
Moreover, there is some $R$ such that each $S_j$ supports no more than $R$ mapping classes. Thus,
\begin{align*}
    \sum_{i=1}^k d_{w_{i-1} \cdot \alpha_i}(\mu, w \cdot \mu) &\leq \sum_{i=1}^k d_{w_{i-1} \cdot S_i}(\mu, w \cdot \mu) + 24\delta \\
    &\leq R(\sum_{S' \subseteq S} [[d_{S'}(\mu, w \cdot \mu)]]_K) \\
    &\leq R(A \cdot d_{\widetilde{\mathcal{M}}(S)}(\mu, w \cdot \mu) + B).
\end{align*}
where the last inequality follows from the Masur-Minsky distance formula.
\end{proof}